\documentclass[11pt]{amsart}

 \usepackage[T1]{fontenc}
\usepackage{enumerate}
\usepackage{geometry}
\usepackage{cite}
\usepackage[colorlinks=true,linkcolor=red,citecolor=blue,hypertexnames=false]{hyperref}

\usepackage{tikz}
\usepackage{graphicx}
\usepackage{pgfplots}
\usepackage{subcaption}
\usepackage{stmaryrd}
\usepackage{microtype}
\usepackage{enumerate}

\usepackage{subcaption}

\definecolor{BC}{HTML}{0000aa}
\definecolor{GV}{HTML}{aa0000}

\usepackage{amsmath}
\usepackage{amsfonts}
\usepackage{amssymb}
\usepackage{amsthm}
\usepackage{mathtools}

\makeatletter
\@namedef{subjclassname@2020}{%
	\textup{2020} Mathematics Subject Classification}
\makeatother

\makeatletter
\newcommand{\customlabel}[2]{%
   \protected@write \@auxout {}{\string \newlabel {#1}{{#2}{\thepage}{#2}{#1}{}} }%
   \hypertarget{#1}{#2}%
}
\newcommand{\customlabelquiet}[2]{%
   \protected@write \@auxout {}{\string \newlabel {#1}{{#2}{\thepage}{#2}{#1}{}} }%
   \hypertarget{#1}{}%
}
\makeatother

\usepackage{fixltx2e}

\newtheorem{theorem}{Theorem}[section]
\newtheorem*{theorem*}{Theorem}
\newtheorem{proposition}[theorem]{Proposition}
\newtheorem*{proposition*}{Proposition}
\newtheorem{corollary}[theorem]{Corollary}
\newtheorem*{corollary*}{Corollary}
\newtheorem{lemma}[theorem]{Lemma}
\newtheorem*{lemma*}{Lemma}
\theoremstyle{definition}
\newtheorem{definition}[theorem]{Definition}
\newtheorem*{definition*}{Definition}
\newtheorem{remark}[theorem]{Remark}
\newtheorem*{remark*}{Remark}

\usepackage{cleveref}

\newcommand{\Ho}{\textrm{\bf (H1)}}
\newcommand{\Ht}{\textrm{\bf (H2)}}
\newcommand{\Htpl}[1]{\textrm{\bf (H2-$#1$)}}

\renewcommand\epsilon{\varepsilon}
\newcommand{\llparen}{(\hspace{-0.55ex}(}
\newcommand{\rrparen}{)\hspace{-0.55ex})}
\newcommand\F{\mathbb{F}}

\newcommand\Q{\mathbb{Q}}
\newcommand\Z{\mathbb{Z}}
\newcommand\N{\mathbb{N}}
\newcommand\PP{\mathbb{P}}
\newcommand\OO{\mathcal{O}}
\newcommand\CC{\mathcal{C}}
\newcommand\HH{\mathcal{H}}
\newcommand\supp{\mathrm{supp}}
\newcommand\Div{\mathrm{Div}}
\newcommand\End{\mathrm{End}}
\newcommand\Hom{\mathrm{Hom}}

\newcommand\pl{\mathfrak{p}} 
\newcommand\ql{\mathfrak{q}} 
\newcommand\rk{\mathrm{rk}}
\newcommand\srk{\mathrm{srk}}
\newcommand\ie{\textit{i.e.~}}

\newcommand\Nrd{N_{\text{\rm rd}}}

\begin{document}
\sloppy
\title{Algebraic Geometry codes in the sum--rank metric}

\author{Elena Berardini}
\address{CNRS; IMB, Université de Bordeaux, 351 cours de la Libération, 33405 Talence, France\\and\\Eindhoven University of Technology, The Netherlands}
\email{elena.berardini@math.u-bordeaux.fr}

\author{Xavier Caruso}
\address{CNRS; IMB, Université de Bordeaux, 351 cours de la Libération, 33405 Talence, France}
\email{xavier.caruso@normalesup.org}


	\maketitle
\begin{abstract}
We introduce the first geometric construction of codes in the sum-rank metric, which we called linearized Algebraic Geometry codes, using quotients of the ring of Ore polynomials with coefficients in the function field of an algebraic curve. We study the parameters of these codes and give lower bounds for their dimension and minimum distance. Our codes exhibit quite good parameters, respecting a similar bound to Goppa’s bound for Algebraic Geometry codes in the Hamming metric. Furthermore, our construction yields codes asymptotically better than the sum-rank version of the Gilbert--Varshamov bound.
\end{abstract}

\setcounter{tocdepth}{1}
\tableofcontents

\section*{Introduction}
Linear codes in the Hamming metric have been playing a central role in 
the theory of error correction since the 50's. Codes in the rank metric, 
firstly introduced by Delsarte for combinatorial interest 
\cite{Delsarte78}, were rediscovered in the last 20 years in the context 
of network coding and, in general, of error correction. Codes in the 
sum-rank metric were introduced more recently. They can be defined as 
follows. Let $k$ be a field. For an integer $s$, let 
$\underline{W}=(W_1,\dots,W_s)$ and $\underline{V}=(V_1,\dots,V_s)$ be two $s$-uples of $k$-vector spaces. 
Write $m_i = \dim_k W_i$ and  $n_i = \dim_k V_i$. Let $\Hom_k(W_i,V_i)$ 
denote the space of $k$-linear morphisms from $W_i$ to $V_i$. We set
\[\Hom_k(\underline{W},\underline{V})\coloneqq \Hom_k(W_1,V_1)\times \dots\times  \Hom_k(W_s,V_s).\]
This is a vector space over $k$ of dimension $\sum_{i=1}^s m_i n_i$.

\begin{definition*}
Let $\underline{\varphi}=(\varphi_1,\dots,\varphi_s)\in \Hom_k(\underline{W}, \underline{V})$. The sum-rank weight of $\underline{\varphi}$ is defined as \[w_{\srk}(\underline{\varphi})\coloneqq\sum_{i=1}^s \rk (\varphi_i)=\sum_{i=1}^s \dim_k \varphi_i(W_i).\]
The sum-rank distance between $\varphi,\psi\in \Hom_k(\underline{W},\underline{V})$ is 

\[d_{\srk}(\underline{\varphi},\underline{\psi})\coloneqq w_{\srk}(\underline{\varphi}-\underline{\psi}). \]
\end{definition*}

 \begin{definition*}
A code $\mathcal{C}$ in the sum-rank metric is a $k$--linear subspace of 
$\Hom_k(\underline{W},\underline{V})$ endowed with the sum-rank distance.
By defintion, its \emph{length} $n$ is $\sum_{i=1}^s m_i n_i$. Its 
\emph{dimension} $\kappa$ is $\dim_k \mathcal{C}$. Its \emph{minimum 
distance} is 
\[d\coloneqq\min \left\{w_{\srk}(\underline{\varphi}) \mid 
\underline{\varphi} \in \mathcal{C}, \underline{\varphi}\neq 
\underline{0}\right\}.\]
\end{definition*}

\begin{remark*} For sake of clarity, we mention that the more standard definition of codes in the sum-rank metric is given in terms of matrices. Let $M_{m,n} (k)$ denote the space of $m \times n$ matrices with coefficients in the field $k$. It is clear that, after the choice of $k$-bases of $W_i$ and $V_i$, we have
\[\Hom_k(\underline{W},\underline{V})=\prod_{i=1}^s M_{m_i,n_i}(k).\]
The definitions of the sum-rank weight and of the sum-rank distance previously introduced can be applied straightforward to elements in $\prod_{i=1}^s M_{m_i,n_i}(k)$. Then, a code in the sum-rank metric is a $k$--linear subspace of $ \prod_{i=1}^s M_{m_i,n_i}(k)$, endowed with the sum-rank distance.
It is clear that the sum-rank distance stays unchanged when considering an element in $ \prod_{i=1}^s M_{m_i,n_i}(k)$ or in $\Hom_k(\underline{W},\underline{V})$, since the rank is invariant under base change.
\end{remark*}

When $n_i = m_i =1$ for all $i\in\{1,\dots,s\}$, the previous definition 
reduces to codes of length $s$ with the Hamming metric and, in the case 
where $s = 1$, to rank-metric codes. However, the sum-rank metric does 
not reduce to a mere generalization of the two aforementioned metrics. 
For instance, codes in the sum-rank metric offer a solution to problems 
in multi-shot linear network coding, space-time coding, and distributed 
storage. We refer the reader to \cite{MSK22} for a detailed introduction 
to the theory of sum-rank metric codes and their applications.

An important case of interest occurs when we are given a finite extension $\ell$ of 
$k$ of degree~$r$, and we set $V_i=\ell$ for every $i$. In this case, $n_i = 
r$ for all~$i$ and the ambient space 
$\Hom_k(\underline{W},\ell)$ is itself a vector space over~$\ell$.
We are then more particularly interested in $\ell$--linear codes which are,
by definition, $\ell$--linear subspaces $\CC \subset \Hom_k(\underline{W}, 
\ell)$. 
Consequently, we define $\ell$--variants of the parameters: the $\ell$--length of $\CC$ is 
$n_\ell \coloneqq
\sum_{i=1}^s m_i$, the $\ell$--dimension of $\CC$ is $\kappa_\ell \coloneqq
\dim_\ell \CC$, and the 
minimal distance $d$ of $\CC$ stays unchanged.
Those three main parameters are related by the equivalent of the 
Singleton bound in the Hamming metric, that in the aforementioned 
setting reads $d + \kappa_\ell \leq n_\ell + 1$ \cite[Prop.~34]{MP18}. Codes
with parameters attaining this bound are called \emph{Maximum Sum-Rank 
Distance (MSRD)}.

\smallskip

Among the most used families of linear codes in the Hamming metric are the Reed--Solomon codes \cite{RS60}. They have parameters attaining the Singleton bound, and benefit from efficient decoding algorithms. The counterpart of Reed--Solomon codes in the rank metric are Gabidulin codes \cite{Gabidulin85}, and linearized Reed--Solomon codes in the sum-rank metric \cite{MP18}. These codes have parameters attaining the analogue of the Singleton Bound in their respective metric, and benefit from efficient decoding algorithms derived from the ones for Reed--Solomon codes \cite{L06,PWZ16,CD18}. 

The main drawback of Reed--Solomon codes is that the storage size of the coordinates of the vectors increases logarithmically with the number of coordinates: in order to have \emph{long} Reed--Solomon codes, one must work over \emph{large} finite fields. The so-called Algebraic Geometry (AG) codes, introduced by V. D. Goppa \cite{Goppa82}, generalize Reed--Solomon codes and benefit from similar properties, while being free of this limitation.  AG codes are constructed by evaluating spaces of functions at rational points on algebraic curves. Since we can find a curve with an arbitrarily large number of rational points, the construction proposed by Goppa yields codes that are generally longer than the Reed--Solomon codes, and thus allows to work on smaller finite fields. AG codes became particularly famous when Tsfasman, Vlăduţ, and Zink used them with modular curves to construct codes with parameters asymptotically better than the Gilbert--Varshamov bound \cite{TVZ82}.

\subsection*{Motivations} 

In contrast with the situation of codes in the Hamming metric, only a few constructions of codes are known in the rank and the sum-rank metric. In particular, no geometric construction has been proposed in these two metrics so far. Furthermore, MSRD codes, such as linearized Reed--Solomon codes, suffer from the same limitation as Reed--Solomon codes. Indeed, keeping the same notation as before, and denoting by $q$ the cardinality of $k$, in \cite[Thm.~6.12]{BGLR21} it is shown that if $\CC\subseteq\Hom_k(\underline{W}, \ell)$ is a MSRD code with minimum distance $d\leq r+2$, then $s\leq q+1$ if $r=1$ and $s\leq q$ otherwise. Furthermore, for MSRD codes of $\ell$--dimension $2$, we have a similar bound, that we prove now.

\begin{lemma*}
Let $\CC\subseteq \End_k(\ell)^s$ be a code in the sum-rank metric of $\ell$-dimension $2$ and of minimum distance $rs - 1$. Then, if $r=1$ we have $s\leq q+1$, otherwise we have $s \leq q-1$.
\end{lemma*}
\begin{proof}
We consider a $\ell$-basis of $\CC$, say $(f_1,\dots, f_s),(g_1,\dots, g_s)$ with $f_i$ and $g_i$ $k$-linear endomorphisms of $\ell$.
For any $i\in\{1,\ldots,s\}$ and any $x \in \ell^\times \coloneqq \ell \setminus \{0\}$, we define the following element of $\PP^1(\ell)$:
\[v_{i,x} \coloneqq [ f_i(x) : g_i(x) ].\]
It is easy to check that for any $a\in k^\times$, we have $v_{i,x} = v_{i,ax}$.

Let us prove that for $i,j\in\{1,\dots,s\}$ and $x, y \in \ell^\times$, we have $v_{i,x}\neq v_{j,y}$, unless $i= j$ and $x$ and $y$ are $k$-collinear.
Indeed, suppose that $v_{i,x} = v_{j,y}$. Then, by construction, the vectors $(f_i(x), g_i(x))$ and $(f_j(y), g_j(y))$ are collinear. Thus, there exist $u, v$, not both zero, such that 
\begin{align*}
(uf_i + vg_i)(x) &= 0,\\
(uf_j + vg_j)(y) &= 0.
\end{align*}
If $i\neq j$, then $uf_i + vg_i$ and $uf_j + vg_j$ are both of rank  smaller than $r$, hence $w_{srk}\left((uf_h+vg_h)_{h\in\{1,\dots,s\}}\right)\leq rs-2$, which contradicts the assumption on the minimum distance of $\CC$. Similarly, if $i = j$, but $x$ is not $k$-collinear to $y$, then we deduce that $uf_i + vg_i$ is of rank at most $r-2$, which again contradicts the hypothesis on the minimum distance of $\CC$. In conclusion, we must have $i = j$ and $x \equiv y \pmod k$.
We infer that the number of pairs $(i,x)$ with $x \pmod k \in \ell^\times$ is at most equal to the cardinal of $\PP^1(\ell)$, \ie
\[s  \frac{q^r - 1}{q - 1} \leq q^r + 1,\]
and hence we have
\begin{equation}\label{eq:t}
s \leq (q - 1) \frac{q^r + 1}{q^r - 1}.
\end{equation}
Since $s$ is necessarily an integer, this implies $s \leq q - 1$ when $q^r>2q-1$, which happens as soon as $r>1$. If $r=1$, then Equation \eqref{eq:t} gives $s\leq q+1$.
\end{proof}
\begin{remark*}
A bound similar to the one stated in the previous lemma can be retrieved using  \cite[Thm.~6.12]{BGLR21} on the dual of the MSRD code of dimension $2$. However, this would give a slightly worse bound, that is $s\leq q$ instead of $s\leq q-1$, when $r>1$. Furthermore, our proof makes use of completely different techniques than the ones developed in the aforementioned paper, and we therefore believe it is of interest on itself.
\end{remark*}

\begin{remark*}
The previous bound could be generalised to codes in $\Hom_k(\underline{W},\ell)$, \ie to codes whose components are rectangular matrices. However, for the purpose of this paper, we focus on the bound for codes sitting in  $\End_k(\ell)^s$,  \ie to codes whose components are square matrices, since this is the setting in which we will develop our geometric codes. In the case of non-square matrices, multiple constructions of MSRD codes longer than linearized Reed--Solomon codes were proposed in \cite{MP22}.
\end{remark*}
\subsection*{Our contribution}
In this paper we present the first geometric construction of codes in the sum-rank metric, from algebraic curves, that we call \emph{linearized Algebraic Geometry codes}.

 Gabidulin and linearized Reed--Solomon codes are constructed using so-called linearized polynomials and Ore polynomials, as introduced by Ore in 1933 \cite{O33}. 
 Taking inspiration from  the approach of \cite{MS98}, where the authors propose a construction of AG codes in the Hamming metric using division algebras over the function field of a curve, in this paper we work with algebras obtained as quotient of rings of Ore polynomials. We develop the theory of Riemann--Roch spaces over Ore polynomials rings with coefficients in the function field of a curve, by exploiting the classical theory of divisors and Riemann--Roch spaces on algebraic curves. On the one hand, this allows us to propose an explicit construction of geometric codes in the sum-rank metric from curves. On the other hand, we can exploit our theory to study the parameters of these new codes.
 
 The geometric codes that we propose are in general longer than linearized Reed--Solomon codes, and have parameters that turn out to respect a similar bound to Goppa's bound for AG codes in the Hamming metric.
Furthermore, we shall prove that, over a finite field of the form $\F_{q^2}$ with $q \geq 11$, our construction yields codes with parameters beating the sum-rank analogue of the Gilbert--Varshamov bound.

\subsection*{Organisation of the paper}

Section~\ref{sec:ore} is devoted to the background on the rings of Ore 
polynomials and to the proofs of some results on the algebras obtained 
as quotients of the ring of Ore polynomial. In Section~\ref{sec:curves}, 
after recalling some general notions on algebraic curves and their 
function fields, we present the theory of Riemann--Roch spaces over the 
rings of Ore polynomials with coefficients in the function field of a 
curve, and we prove a bound on their dimension using the classical 
Riemann--Roch theorem. In Section~\ref{sec:thecode}, we construct linearized Algebraic Geometry codes, study their parameters and their
asymptotic behaviour with regard to the Gilbert--Varshamov bound. When 
considering the case of the projective line, we retrieve linearized 
Reed--Solomon codes. Finally, Section~\ref{sec:conclusion} serves as a 
general conclusion in which we compare our results with those 
of~\cite{MS98}, and discuss several perspectives.

\section{Ore polynomial rings}\label{sec:ore}

Throughout this section, we fix three positive integers~$r$, $d$
and $m$ such that $r = md$.
We consider a field $K$, together with a Galois extension $L_0/K$
such that $\text{Gal}(L_0|K) \simeq \Z/d\Z$. We denote by $\Phi_0$ 
a generator of the former Galois group. We set $L \coloneqq L_0^m$, equipped with the coordinate-wise product, and embed
$K$ and $L_0$ diagonally into $L$, that is, sending an element $a$ to $(a,\dots,a)$. We define
$$\Phi : L \to L, \quad (a_1, \ldots, a_m) \mapsto
\big(\Phi_0(a_m), a_1, \ldots, a_{m-1}\big).$$
It is easy to check that $\Phi$ has order~$r$, and that an element
$x \in L$ is a fixed point of $\Phi$ if and only if $x$ lies in~$K$.
Besides, an element $a = (a_1, \ldots, a_m) \in L$ is invertible if
and only if $a_i \neq 0$ for all $i$. We write $L^\times$ for the
subset of invertible elements of $L$.

We denote by $N_{L_0/K} : L_0 \to K$ the norm map of $L_0$ over $K$
and similarly, for $a = (a_1, \ldots, a_m) \in L$, we set
$$N_{L/K}(a) \coloneqq N_{L_0/K}(a_1) \cdot N_{L_0/K}(a_2) \cdots
N_{L_0/K}(a_m).$$
This defines a multiplicative function $N_{L/K} : L \to K$ which
sends $L^\times$ to $K^\times \coloneqq K \setminus \{0\}$.

\medskip
The goal of the rest of the section is to present some results from the standard theory of central simple algebras. A classical reference in this context is \cite{Reiner75}. Here we present the results in our setting, that is for the ring of Ore polynomials, which allows us to give more concise proofs.
\subsection{The algebra $D_{L,x}$}
\label{ssec:DLx}

We let $L[T; \Phi]$ denote the ring of Ore polynomials in the
variable $T$. We recall briefly that elements of $L[T; \Phi]$ are
usual polynomials with usual addition; however, the multiplication
on $L[T; \Phi]$ is twisted by the rule $T \cdot a = \Phi(a) T$ for
all $a \in L$.

Given an element $x \in K$, $x \neq 0$, we define
$$D_{L, x} \coloneqq L[T;\Phi] / (T^r - x).$$
Since $T^r - x$ commutes with every element in $L[T; \Phi]$, the quotient $D_{L,x}$
inherits a ring structure.

\begin{lemma}\label{lem:isomDLx} 
\begin{enumerate}[(i)]
\item\label{isoun} Let $u \in L^\times$ and write $v = N_{L/K}(u)$. 
We have an isomorphism of rings
$$\gamma_u : D_{L,x} \stackrel\sim\longrightarrow D_{L,v^{-1}x},
\quad T \mapsto uT.$$
\item\label{isodeux} When $x = 1$, we have an isomorphism of rings
$$\varepsilon : D_{L,1} \stackrel\sim\longrightarrow \End_K(L),
\quad T \mapsto \Phi.$$
\end{enumerate}
\end{lemma}

\begin{proof}
It is easily checked that $\gamma_u$ is a well-defined ring
homomorphism and that $\gamma_{u^{-1}}$ is its inverse. This
proves~\eqref{isoun}.
For~\eqref{isodeux}, we first notice that $\varepsilon$ is well-defined,
given that $\Phi$ has order~$r$. Injectivity boils down to
proving that the family $\left\{\text{Id}, \Phi, \ldots, \Phi^{r-1}\right\}$
is free over $L$, which is a direct consequence of Artin's
theorem on linear independence of characters. Surjectivity 
follows by comparing dimensions (over $K$).
\end{proof}

The quotient rings $D_{L,x}$ are equipped with a so-called
\emph{reduced norm map} $\Nrd : D_{L,x} \to K$ that we define 
now.
For this, we first observe that $D_{L,x}$ is a free $L$-module of 
rank~$r$ with basis~$(1, T, \ldots, T^{r-1})$.

\begin{definition}
Let $f \in D_{L,x}$. The reduced norm of $f$, denoted by $\Nrd(f)$,
is the determinant of the $L$-linear map $D_{L,x} \to D_{L,x}$,
$g \mapsto gf$.
\end{definition}

Concretely, the reduced norm of $f = a_0 + a_1 T + \ldots
a_{r-1} T^{r-1} \in D_{L,x}$ is the determinant of the matrix
\begin{equation}
\label{eq:matrixMf}
M_f = \left(
\begin{matrix}
a_0 & x{\cdot}\Phi(a_{r-1}) & \cdots & x{\cdot}\Phi^{r-1}(a_1) \\
a_1 & \Phi(a_0) & \cdots & x{\cdot}\Phi^{r-1}(a_2) \\
\vdots & \vdots & & \vdots \\
a_{r-2} & \Phi(a_{r-3}) & \cdots & x{\cdot}\Phi^{r-1}(a_{r-1}) \\
a_{r-1} & \Phi(a_{r-2}) & \cdots & \Phi^{r-1}(a_0)
\end{matrix}
\right).
\end{equation}
One readily checks that
$$\Phi(M_f) = 
\left(\begin{matrix}
& 1 \\
& & \ddots \\
& & & 1 \\
x^{-1}
\end{matrix}\right) \cdot
M_f \cdot
\left(\begin{matrix}
& & & x \\
1 \\
& \ddots \\
& & 1
\end{matrix}\right),$$
from what we deduce that $\Nrd(f) = \det(M_f)$ is invariant under
$\Phi$; hence $\Nrd(f) \in K$ as we claimed earlier.
Moreover, the reduced norm map behaves well with respect to the 
isomorphisms $\gamma_u$ and $\varepsilon$ of Lemma~\ref{lem:isomDLx}, as showed in the following lemma.

\begin{lemma}
\label{lem:NrdDLx}
\begin{enumerate}[(i)]
\item\label{lem:un} For $f \in D_{L,x}$ and $u \in L^\times$, we have
$\Nrd(f) = \Nrd(\gamma_u(f))$.
\item\label{lem:deux}  For $f \in D_{L,1}$, we have
$\Nrd(f) = \det (\varepsilon(f))$.
\end{enumerate}
\end{lemma}

\begin{proof}
\eqref{lem:un}~Write $v = N_{L/K}(u)$ and let $\mu$ (resp $\mu'$) be the $L$-linear 
endomorphism of $D_{L,x}$ (resp. $D_{L,v^{-1}x}$) taking $g$ to $gf$
(resp. to $g{\cdot}\gamma_u(f)$). The maps $\mu$ and $\mu'$ are conjugated
under the isomorphism $\gamma_u$. Hence, their determinants agree, showing
that $\Nrd(f) = \Nrd(\gamma_u(f))$.

\eqref{lem:deux}~For $f \in L$ (resp. $f \in D_{L,1}$), let $\mu_f$ denote the right
multiplication by $f$ on $L$ (resp. on $D_{L,1}$).
We consider the tensor product $L \otimes_K L$ and view it as a 
$L$-vector space by letting $L$ act on the first factor. Since $L/K$ is
Galois with cyclic group generated by $\Phi$, it follows from Galois
theory that we have a $L$-linear decomposition
$$\begin{array}{rcl}
L \otimes_K L & \stackrel\sim\longrightarrow& L^r \\
x \otimes y & \mapsto &\big(x{\cdot} \Phi^i(y)\big)_{0 \leq i < r}.
\end{array}$$
In the corresponding $L$-basis of $L \otimes_K L$,
the matrices of the endomorphisms
$1 \otimes \mu_a$ and $1 \otimes \Phi$ are, respectively,
$$\left(\begin{matrix}
a \\
& \Phi(a) \\
& & \ddots \\
& & & \Phi^{r-1}(a)
\end{matrix}\right)
\quad \text{and} \quad
\left(\begin{matrix}
& & & 1 \\
1 \\
& \ddots \\
& & 1
\end{matrix}\right).$$
Hence, the matrix of the endomorphism $1 \otimes \varepsilon(f)$ is
exactly the matrix $M_f$ defined in 
Equation~\eqref{eq:matrixMf}. The equality $\Nrd(f) = \det (\varepsilon(f))$ 
follows.
\end{proof}

\subsection{Over Laurent series rings}
\label{ssec:DLxlocal}

In this subsection, we fix a field $k$, and write $\OO_K \coloneqq k\llbracket t \rrbracket$ for the ring of power series in $t$. We set $K = \OO_K[1/t]=k\llparen t \rrparen$. We recall that $K$ is the field of fractions of $\OO_K£$, called the \emph{field of Laurent series} over~$k$. 
The elements in $K$ are series in $t$ with possibly a finite number of negative exponents.
The smallest exponent appearing in a Laurent series $f$ is called its \emph{$t$-adic valuation} and is denoted by $v_t(f)$. This defines a function
$v_t : K \to \Z \sqcup \{\infty\}$.
The former extends uniquely to a valuation on $L_0$ that, in a slight abuse
of notation, we continue to denote by $v_t$. Recall that, in full
generality, $v_t$ does not take integral values on $L_0$; more
precisely, if $e$ denotes the ramification index of $L_0/K$, $v_t$
defines a surjective function from $L_0$ to $\frac 1 e \Z \sqcup 
\{\infty\}$. Note that $e$ divides $d$ and hence $r$.
We write $\OO_{L_0}$ for the valuation ring of $L_0$, that is the 
subring of $L_0$ formed by elements with nonnegative valuation.

We recall that we have defined $L = L_0^m$. We set $\OO_L \coloneqq
(\OO_{L_0})^m$ accordingly.
For $j \in \{1, \ldots, m\}$, we consider the function $v_{j,t}$ on
$L$ defined by $v_{j,t}(c_1, \ldots, c_m) \coloneqq v_t(c_j)$ for $c_1,
\ldots, c_m \in L_0$.
Similarly, for $f = a_0 + a_1 T + \cdots + a_{r-1} T^{r-1} \in 
D_{L,x}$ (with $a_i \in L$), we set
$$w_{j,x}(f) \coloneqq \min_{0 \leq i < r}
\left(v_{j,t}(a_i) + i \cdot \frac{v_t(x)} r\right).$$
This defines a function $w_{j,x} : D_{L,x} \to \frac 1 r \Z \sqcup
\{\infty\}$ (for $1 \leq j \leq m$). We further define
$w_x \coloneqq \min_{1 \leq j \leq m} w_{j,x}$. 
Given $f, g \in D_{L,x}$, one checks that:
\begin{itemize}
\item $w_{j,x}(f+g) \geq \min\left(w_{j,x}(f), w_{j,x}(g)\right)$
for $1 \leq j \leq m$,
\item $w_x(f+g) \geq \min\left(w_x(f), w_x(g)\right)$,
\item $w_x(fg) \geq w_x(f) + w_x(g)$,
\item $w_x(f) = \infty$ if and only if $f = 0$.
\end{itemize}
We define $\Lambda_{L,x}$ as the subset of $D_{L,x}$ consisting
of elements $f$ for which $w_{j,x}(f) \geq 0$ for all~$j$; it is
a subring of $D_{L,x}$.

\begin{lemma}
\label{lem:isomLambdaLx}
Let $u = (u_1, \ldots, u_m) \in L^\times$ and set
$y = N_{L/K}(u)^{-1} {\cdot} x$.
Let $\gamma_u : D_{L,x} \to D_{L,y}$ be the isomorphism defined
in Lemma~\ref{lem:isomDLx}.\eqref{isoun}.
If $v_t(u_1) = \cdots = v_t(u_m)$, then
$$w_{j,x}(f) = w_{j,y}(\gamma_u(f))$$
for all $j \in \{1, \ldots, m\}$ and all $f \in D_{L,x}$.
In particular, $\gamma_u$ induces an isomorphism
$\Lambda_{L,x} \stackrel\sim\longrightarrow \Lambda_{L,y}$.
\end{lemma}

\begin{proof}
For simplicity, write $v$ for the common value of $v_t(u_1), \ldots,
v_t(u_m)$. The relation $y = N_{L/K}(u)^{-1} {\cdot} x$ then implies
that $v_t(y) = v_t(x) - r{\cdot}v$.

Consider now an element 
$f = a_0 + a_1 T + \cdots + a_{r-1} T^{r-1} \in D_{L,x}$.
By definition,
$$\gamma_u(f) = 
  \sum_{i=0}^{r-1} a_i \:u \: \Phi(u) \cdots \Phi^{i-1}(u) \cdot T^i.$$
For all $i$ and $j$, we have 
$v_{j,t}(a_i \:u \: \Phi(u) \cdots \Phi^{i-1}(u)) = v_{j,t}(a_i) + 
i{\cdot}v$. Hence
\begin{align*}
w_{j,y}(\gamma_u(f)) 
& = \min_{0 \leq i < r} 
    \left(v_{j,t}(a_i) + i{\cdot}v + i\cdot\frac{v_t(y)}r\right) \\
& = \min_{0 \leq i < r} 
    \left(v_{j,t}(a_i) + i\cdot\frac{v_t(x)}r\right) = w_{j,x}(f),
\end{align*}
which proves the lemma.
\end{proof}

\begin{proposition}
\label{prop:divisionlocal}
We assume that $m=1$, $L_0/K$ unramified and $\gcd(v_t(x), r) = 1$.
Then, $D_{L,x}$ has no nonzero zero divisor.
\end{proposition}

\begin{proof}
Let $f$ and $g$ be nonzero elements in $D_{L,x}$. We want
to prove that $fg$ cannot vanish.
By our assumption on $m$, we have $v_{1,t} = v_t$ and $w_{1,x} = w_x$.
We claim that the minimum in the definition of $w_x(f)$ is reached 
only once; in other words, if $f$ is written as $f = a_0  + a_1 T+ \cdots +
a_{r-1} T^{r-1}$ (with $a_i \in L$), there exists a unique index
$i_f \in \{0, \ldots, r{-}1\}$ such that
$$v_t(a_{i_f}) + i_f \cdot \frac{v_t(x)} r = w_x(f).$$
Indeed, given that $v_t(a_i)$ is an integer for all $i$ by the 
assumption on the ramification, such an index $i_f$ has to satisfy
the congruence $i_f {\cdot} v_t(x) \equiv r {\cdot} w_x(f) \pmod r$.
The latter has a unique solution, given that $v_t(x)$ is coprime with 
$r$. As a conclusion, we can write $f = c_f T^{i_f} + f_1$
where $c_f \in L$ and $f_1 \in D_{L,x}$ satisfy
$w_x(c_f T^{i_f}) = w_x(f)$ and $w_x(f_1) > w_x(f)$.

Similarly, $g = c_g T^{i_g} + g_1$
where $i_g$ is an integer in the range $[0,r)$, and $c_g \in L$ and 
$g_1 \in D_{L,x}$ are such that $w_x(c_g T^{i_g}) = w_x(g)$ and 
$w_x(g_1) > w_x(g)$.
Computing the product $fg$, we find
$$fg = c_f \: \Phi^{i_f}(c_g) T^{i_f + i_g} + h_1,$$
with $w_x(h_1) > w_x(f) + w_x(g)$. On the other hand, we have 
\begin{align*}
w_x\big(c_f \: \Phi^{i_f}(c_g) T^{i_f + i_g}\big) 
& = v_t\big(c_f \: \Phi^{i_f}(c_g)\big) + (i_f + i_g) \cdot \frac{v_t(x)} r \\
& = v_t(c_f) + v_t(c_g) + (i_f + i_g) \cdot \frac{v_t(x)} r\\
  &= w_x(f) + w_x(g).
\end{align*}
Therefore $c_f \: \Phi^{i_f}(c_g) T^{i_f + i_g}$ cannot be equal to
$-h_1$ (because the valuations differ), showing eventually that 
$fg \neq 0$, as wanted.
\end{proof}

We now examine the relationships between the valuations and the
reduced norm.

\begin{proposition}
\label{prop:vtNrd1}
For all $f \in D_{L,x}$, we have
\[\displaystyle v_t\big(\Nrd(f)\big) \geq d \cdot \sum_{j=1}^m w_{j,x}(f).\]
\end{proposition}

\begin{proof}
Write $f = a_0 + a_1 T + \cdots + a_{r-1} T^{r-1}$ with $a_i \in L$.
Let $M_f$ be the matrix defined by Equation~\eqref{eq:matrixMf}, and, for
$1 \leq u, v \leq r$, let $m_{u,v}$ denote its entry in position $(u,v)$.
By definition,
$$m_{u,v} =\begin{cases}
  \Phi^{v-1}(a_{u-v}) & \text{if } u \geq v, \smallskip \\
  x \cdot \Phi^{v-1}(a_{u-v+r}) & \text{otherwise.}
  \end{cases}$$
We extend the valuation $v_t$ on $L$ by
$$v_t(a) = \frac 1 m \sum_{j=1}^m v_{j,t}(a)
\qquad (a \in L).$$
We observe that $v_t$ agrees on $L_0$ (embedded diagonally in $L$)
with the valuation $v_t$ we have defined previously; hence no risk
of confusion is possible. Additionally, $v_t$ is invariant under 
$\Phi$. Note, however, that $v_t$ is no longer a group morphism but
it only satisfies $v_t(ab) \geq v_t(a) + v_t(b)$ for $a, b \in L$.

From the previous properties, we derive
$$v_t(m_{u,v}) =\begin{cases}
  v_t(a_{u-v}) & \text{if } u \geq v, \smallskip \\
   v_t(a_{u-v+r}) + v_t(x) & \text{otherwise.} \\
   \end{cases}$$
On the other hand, it follows from the definition of $w_{j,x}$
that $v_{j,t}(a_i) \geq w_{j,x}(f) - i {\cdot} \frac{v_t(x)} r$
for all $i$ and $j$. Summing over $j$, we get
$$v_t(a_i) \geq
\frac 1 m \sum_{j=1}^{m} w_{j,x}(f)  - i {\cdot} \frac{v_t(x)} r ,$$
and finally
$$v_t(m_{u,v}) \geq
\frac 1 m \sum_{j=1}^{m} w_{j,x}(f) +(v{-}u) {\cdot} \frac{v_t(x)} r $$
for all $u, v \in \{1, \ldots, r\}$.
If $\sigma$ is a permutation of $\{1, \ldots, r\}$, we then find
$$v_t\left(\prod_{u=1}^r m_{u,\sigma(u)}\right) \geq
d \cdot \sum_{j=1}^{m} w_{j,x}(f) +
\frac{v_t(x)} r \sum_{u=1}^r (\sigma(u) - u).$$
The last sum vanishes given that $\sigma$ is a permutation.
We conclude that
$$v_t\big(\Nrd(f)\big) =
v_t\big(\!\det (M_f)\big) \geq d \cdot \sum_{j=1}^{m} w_{j,x}(f),$$
and the proposition is proved.
\end{proof}

To conclude this subsection, we focus on the case $x = 1$.
We recall from Lemma~\ref{lem:isomDLx} that we have an isomorphism
$\varepsilon : D_{L,1} \to \End_K(L)$, defined by $T \mapsto
\Phi$. It follows from the definitions that, when $f \in
\Lambda_{L,1}$, the map $\varepsilon(f)$ takes $\OO_L$ to
itself. Hence $\varepsilon$ induces a ring homomorphism
$\Lambda_{L,1} \to \End_{\OO_K}(\OO_L)$. Reducing modulo $t$
on the right hand side, we get a third ring homomorphism
$\bar \varepsilon : \Lambda_{L,1} \to \End_k(\OO_L/t\OO_L)$.

\begin{proposition}
\label{prop:vtNrd2}
For all $f \in \Lambda_{L,1}$, we have
\[v_t\big(\Nrd(f)\big) \geq \dim_k \ker \bar\varepsilon(f).\]
\end{proposition}

\begin{proof}
Write $\kappa = \dim_k \ker \bar\varepsilon(f)$, and let
$\bar e_1, \ldots, \bar e_r$ be a $k$-basis of $\OO_L/t\OO_L$ 
such that $\bar \varepsilon(\bar e_i) = 0$ for $i \in \{1,
\ldots, \kappa\}$. For $1 \leq i \leq r$, we choose an element
$e_i \in \OO_L$ whose reduction modulo $t$ is $\bar e_i$.
By Nakayama's lemma \cite[Cor.~4.8]{Eisenbud13}, the family $(e_1, \ldots, e_r)$ is a basis
of $\OO_L$ over $\OO_K$. Let $M$ be the matrix of $\epsilon(f)$
in this basis. It has all entries in $\OO_K$, while the first 
$\kappa$ columns of $M$ have entries divisible by~$t$ by construction.
Therefore, $\det(M)$ is divisible by $t^\kappa$.
Finally, we know from Lemma~\ref{lem:NrdDLx} that $\Nrd(f) =
\det(M)$. The proposition follows.
\end{proof}

  \section{Algebraic curves}\label{sec:curves}

  Throughout this section, we let $k$ be a field.
  \subsection{Divisors on curves and Riemann--Roch spaces}

  In this subsection, we recall some classical definitions and results on algebraic curves, and refer the reader to \cite{S09} for a nice  exposition of this theory.

  We consider a smooth projective irreducible algebraic curve $X$ of genus $g_X$ defined over $k$ and we set $K=k(X)$ to be its function field.
  We denote by $X^\star$ the set of places (or, equivalently, closed points) of $X$. 
  Given $\pl \in X^\star$, we let $K_\pl$ be the completion of $K$ at the place $\pl$~\cite[\S 4.2]{S09}.
  It is equipped with the $\pl$-adic valuation $v_\pl$. We denote by $\mathcal{O}_\pl$ its ring of integers and by
  $k_\pl = \OO_\pl/\pl$ its residue class field. 
  The \emph{degree} of $\pl$, denoted by $\deg_X(\pl)$ in what follows, is by definition the degree of the extension
  $k_\pl/k$.
  
   \begin{definition}
 The divisor group of $X$, $\Div(X)$, is the free abelian group generated by the places of $X$. A \emph{divisor} on $X$ is therefore a formal sum 
  \[D=\sum_{\pl\in X^\star} n_\pl \pl \quad\text{ with }n_\pl\in\Z \text{ almost all zero}.\]
 The \emph{degree} of $D$ is defined by $\deg_X(D)=\sum_{\pl\in X^\star} n_\pl \deg_X (\pl)$ and its \emph{support} is $\supp(E)=\{\pl\in X^\star\;|\; n_\pl\neq 0\}$. The divisor $E$ is called \emph{effective}, written $E\geq 0$, if $n_\pl \geq 0$ for any $\pl$. Two divisors are added coefficientwise.
 \end{definition}
 
 The \emph{principal divisor} associated to a rational nonzero function $x\in K$ is \[(x)=\sum_{\pl\in X^\star} v_\pl(x)\:\pl.\]
 Since rational functions in $K$ have the same number of zeros and poles, counted with multiplicity, principal divisors have zero degree.
 
For a divisor $D\in \Div(X)$, we define the associated Riemann--Roch space as
\begin{equation}\label{eq:RRspace}
L_X(D)\coloneqq \{ x\in K^\times \mid (x)+D\geq 0\}\cup\{0\}.
\end{equation}
This is a $k$-vector space of finite dimension.
\begin{theorem}[Riemann--Roch theorem]\label{th:classicalRR}
For any divisor $D\in\Div(X)$ we have
\[\dim_k L_X(D)=\deg_X(D) +1-g_X +\dim_k L_X(K_X-D),\]
where $K_X$ denotes a canonical divisor on $X$.
\end{theorem}
\medskip

\subsection{Riemann--Roch spaces in Ore polynomial rings}
\label{ssec:RROre}

In what follows, we will work extensively with Galois coverings of 
curves and their corresponding extensions at the level of function
fields. For definitions and classical results on this subject, we refer 
the reader to \cite[Chapter~3]{S09}.

We consider two smooth projective irreducible algebraic curves $X$ 
and $Y$ defined over $k$, together with a surjective map $\pi : 
Y \to X$. 
Let $K \coloneqq k(X)$ and $L \coloneqq k(Y)$ denote the fields of functions of $X$
and $Y$ respectively. The map $\pi$ induces a ring homomorphism $K 
\to L$, turning $L$ into an extension of $K$. 
Moreover, we assume that $L/K$ is Galois with
cyclic Galois group of order~$r$. 

\medskip

We denote by $\Div(X)$ and $\Div(Y)$ the group of divisors on $X$ 
and $Y$, respectively, and we set $\Div_\Q(Y) \coloneqq \Div(Y) \otimes \Q$.
To avoid confusion, we reserve the letter $\pl$ (resp. $\ql$) to
denote places of $X$ (resp. of $Y$). We say that a place $\ql$
\emph{divides} $\pl$ or, equivalently, that $\ql$ is \emph{above}
$\pl$, and we note $\ql | \pl$, when $\pi$ maps $\ql$ to $\pl$.
Let $\pl \in X^\star$ and let $K_\pl$ be the completion of $K$ at
the place $\pl$. We have the decomposition
\begin{equation}
\label{eq:decompLp}
K_\pl \otimes_K L \simeq \prod_{\ql | \pl} L_\ql.
\end{equation}
For simplicity, we write $L_\pl = K_\pl \otimes_K L$.

We fix a generator $\Phi \in \text{Gal}(L|K)$. For any place
$\pl \in X^\star$, we note that $\Phi$ permutes cyclically the
$L_\ql$'s of Equation~\eqref{eq:decompLp}. Hence, they are all
isomorphic and one can number the places above $\pl$ as follows
$$\pi^{-1}(\pl) = \{\ql_1, \ql_2, \ldots, \ql_{m_\pl}\},$$
in such a way that $\Phi$ maps $L_{\ql_j}$ to $L_{\ql_{j+1}}$ for 
all $j$ (with the convention that $\ql_{m_\pl+1} = \ql_1$).
The morphism $\Phi^{m_\pl}$ then induces an automorphism $\Phi_{\pl, 0}$
of $L_{\ql_1}$ of order~$d_\pl = r/m_\pl$.
Setting $L_{\pl, 0} = L_{\ql_1}$, we finally see that the pair
$(L_\pl, K_\pl)$ fits in the framework of Subsection~\ref{ssec:DLxlocal}.
\medskip

Let $x$ be a fixed function in $K^\times$.
We consider the algebras $D_{L,x} = L[T; \Phi]/(T^r - x)$ and
$D_{L_\pl,x} = L_\pl[T; \Phi]/(T^r - x)$. We recall that we have
defined in Subsection~\ref{ssec:DLxlocal} the valuations
$$\textstyle
w_{j,x} : D_{L_\pl,x} \to \frac 1 r \Z \sqcup \{\infty\}
\quad (1 \leq j \leq m_\pl).$$
Instead of indexing them by the integers $j \in \{1, \ldots, m_\pl\}$,
it is more convenient here to index them by the places above $\pl$, 
\ie writing $w_{\ql_j, x}$ for $w_{j,x}$.
For an element $f \in D_{L_\pl, x}$ written as 
$f = f_0 + f_1 T + \cdots + f_{r-1} T^{r-1}$ (with $f_i \in L_\pl$),
we then have
$$w_{\ql, x}(f) =
\min_{0 \leq i < r} \left(\frac{v_\ql(f_i)}{e_\ql} + i \cdot \frac{v_\pl(x)} r
\right),$$
where $e_\ql$ denotes the ramification index at $\ql$, which is
also the ramification index of the extension $L_\ql/K_\pl$. Since
all the $L_\ql$'s are isomorphic, we see that $e_\ql$ depends only on the
place $\pl$ below; for this reason, we will often denote it by $e_\pl$
in what follows.

For a place $\pl \in X^\star$, we set
$$\rho_\pl = \frac{e_\pl \cdot v_\pl(x)} r,$$
and define $a_\pl$ and $b_\pl$ by $\rho_\pl = \frac{a_\pl}{b_\pl}$, 
where the latter fraction is irreducible and its denominator $b_\pl$
is positive. Since $v_\pl(x)$ vanishes for almost all places $\pl$, we 
find that $\rho_\pl = 0$, $a_\pl = 0$ and $b_\pl = 1$ for almost all 
$\pl \in X^\star$.

\begin{definition}[Riemann--Roch spaces of $D_{L,x}$]
\label{def:RROre}
Let
$E = \sum_{\ql \in Y^\star} n_\ql \ql \in \Div_\Q(Y)$
where, for all $\ql$, the coefficient $n_\ql$ is in $\frac 1{b_\pl}
\Z$ where $\pl = \pi(\ql)$ is the place below $\ql$. We define the \emph{Riemann--Roch space} of $D_{L,x}$ associated with $E$ as
$$\Lambda_{L,x}(E) \coloneqq \big\{\, f \in D_{L,x} \,|\,
e_\ql w_{\ql,x}(f) + n_\ql \geq 0 \text{ for all } \ql \in Y^\star \,\big\}.$$
\end{definition}

\begin{remark}
We use the letter $E$ (instead of $D$) to denote the divisor in 
order to lower the risk to create confusion with the algebra $D_{L,x}$.
\end{remark}

Keeping the notation of Definition~\ref{def:RROre}, it follows readily 
from the definitions that
\begin{equation}
\label{eq:decompRROre}
\Lambda_{L,x}(E) = \bigoplus_{i=0}^{r-1} L_Y(E_i) \cdot T^i,
\end{equation}
where, letting $\lfloor \cdot \rfloor$ denote the lower integer part function,
the divisors $E_i$ are defined by
$$E_i \coloneqq \sum_{\ql \in Y^\star} \big\lfloor
n_\ql + i \cdot \rho_{\pi(\ql)} \big\rfloor \cdot \ql
\in \Div(Y) \qquad (0 \leq i < r),$$
and the $L_Y(E_i)$'s are the ``classical'' Riemann--Roch spaces (on $Y$), as defined in Equation~\eqref{eq:RRspace}.

\begin{lemma}
\label{lem:degEi}
We have
\[\displaystyle \sum_{i=0}^{r-1} \deg_Y(E_i) = r {\cdot} \deg_Y(E)
- \frac {r^2} 2 \sum_{\pl \in X^\star} \frac{b_\pl{-}1}{b_\pl e_\pl} \deg_X(\pl).\]
\end{lemma}

\begin{proof}
Fix a place $\ql \in Y^\star$, and write $\pl = \pi(\ql)$ and
$n_\ql = \frac{c_\ql}{b_\pl}$. For $i \in \{0, \ldots, r{-}1\}$,
we have
$$\big\lfloor n_\ql + i \cdot \rho_{\pi(\ql)} \big\rfloor
= \left\lfloor \frac{c_\ql + i {\cdot} a_\pl}{b_\pl} \right\rfloor
= \frac{c_\ql + i {\cdot} a_\pl - \varepsilon_{i,\ql}}{b_\pl},$$
where $\varepsilon_{i,\ql}$ denotes the remainder in the division 
of $c_\ql + i {\cdot} a_\pl$ by $b_\pl$. 
From the fact that $a_\pl$ and $b_\pl$ are coprime, we derive
that for each value $\varepsilon \in \{0, \ldots, b_\pl{-}1\}$,
there are exactly $\frac r{b_\pl}$ indices $i$ for which
$\varepsilon_{i,\ql} = \varepsilon$.
Therefore, summing over $i$, we get
\begin{align*}
\sum_{i=0}^{r-1}
\big\lfloor n_\ql + i \cdot \rho_{\pi(\ql)} \big\rfloor
& = r {\cdot} n_\ql + \frac{r(r{-}1)}2 {\cdot} \rho_{\pi(\ql)} - 
  \frac{r(b_\pl{-}1)}{2b_\pl} \\
& = r {\cdot} n_\ql + \frac{r{-}1}2 {\cdot} v_\ql(x) - 
  \frac{r(b_\pl{-}1)}{2b_\pl}.
\end{align*}
Summing over $\ql$ and weighting by $\deg_Y(\ql)$, and using that $\sum_{\ql\in Y^\star} v_\ql(x)\deg_Y(\ql)=0$, we end up with
$$\sum_{i=0}^{r-1} \deg_Y(E_i) = r \cdot \deg_Y(E) - 
\frac r 2 \sum_{\ql \in Y^\star} \frac{b_{\pi(\ql)} -1}{b_{\pi(\ql)}}\deg_Y(\ql).$$
Noticing finally that $\deg_Y(\ql) = \frac r{m_\pl e_\pl} {\cdot} \deg_X(\pi(\ql))$,
we obtain the announced formula.
\end{proof}

\begin{corollary}
\label{cor:RROre}
For a divisor $E = \sum_{\ql \in Y^\star} n_\ql \ql \in \Div_\Q(Y)$ as
in Definition~\ref{def:RROre}, the space $\Lambda_{L,x}(E)$ is finite
dimensional over $k$ and
\begin{align*}\dim_k \Lambda_{L,x}(E) \geq&
r {\cdot} \deg_Y(E) - r{\cdot}(g_Y - 1)+\\
&-\frac {r^2} 2 \sum_{\pl \in X^\star} \frac{b_\pl{-}1}{b_\pl e_\pl} \deg_X(\pl).
\end{align*}
\end{corollary}

\begin{proof}
On the one hand, from Equation~\eqref{eq:decompRROre}, we derive
$$\dim_k \Lambda_{L,x}(E) = \sum_{i=0}^{r-1} \dim_k L_Y(E_i).$$
On the other hand, it follows from the classical Riemann--Roch theorem (Theorem \ref{th:classicalRR})
that $\dim_k L_Y(E_i) \geq \deg_Y E_i - (g_Y - 1)$. Combining this input
with Lemma~\ref{lem:degEi}, we get the corollary.
\end{proof}

\begin{remark}
We point out that equality in the bound of Corollary \ref{cor:RROre} is attained whenever for any $i$ we have $\dim_k L_Y(E_i) = \deg_Y E_i - (g_Y - 1)$, which happens as soon as $\deg_Y E_i \geq 2g_Y -1$ for any $i$. 
\end{remark}

\section{Linearized Algebraic Geometry codes}\label{sec:thecode}

In this section we introduce codes in the sum-rank metric from algebraic 
curves, that we call linearized Algebraic Geometry codes.  We propose a 
general construction using a morphism $\pi : Y \to X$ between two curves 
(Subsection \ref{ssec:construction}). We give bounds for the dimension 
and the minimum distance of our codes (Theorem \ref{th:codeparameters}). 
In Subsection~\ref{ssec:isotrivial}, we consider the case of isotrivial 
covers. In particular, when the curve $X$ has genus $g_X=0$, we retrieve 
the construction of linearized Reed--Solomon codes, as proposed in 
\cite{MP18, CD22}. Finally, in Subsection~\ref{ss:GV}, we
study the asymptotic behaviour of our codes in the isotrivial case
and compare it with the Gilbert--Varshamov bound.

\subsection{The code construction}\label{ssec:construction}

We consider the setting of Subsection~\ref{ssec:RROre} and keep all the
notation from here. In particular, we fix a base field~$k$ and consider
a map $\pi : Y \to X$ between 
smooth projective irreducible algebraic curves defined over $k$. 
We write $K \coloneqq X(k)$ and $L \coloneqq Y(k)$ for the function fields of $X$ and $Y$, respectively. We assume that the
extension $L/K$ is Galois with cyclic Galois group of order~$r$,
generated by $\Phi$. As in Subsection~\ref{ssec:RROre}, we continue
to use the letter $\pl$ (resp. $\ql$) to refer to places of $X$
(resp. of $Y$). We fix in addition:
\begin{itemize}
\item a function $x \in K^\times$,
\item a divisor $E = \sum_{\ql \in Y^\star} n_\ql \ql \in \Div_\Q(Y)$
satisfying the condition of Definition~\ref{def:RROre},
\item a positive integer $s$ and $s$ \emph{rational} places
$\pl_1, \ldots, \pl_s \in X^\star$ which do not belong to
$\pi(\supp(E))$.
\end{itemize}
For $i \in \{1, \ldots, s\}$, we write $K_i \coloneqq K_{\pl_i}$ (the
completion of $K$ at the place $\pl_i$) and set $L_i \coloneqq K_i \otimes_K 
L$. Since the $\pl_i$'s are rational, we have an isomorphism
$K_i \simeq k\llparen t_i \rrparen$, where $t_i$ is a uniformizing
parameter at $\pl_i$. We let $m_i$ be the number of places above
$\pl_i$.

\medskip

We formulate several hypotheses, the second one depending on a
place $\pl$ of $X$:

\begin{itemize}
\item[\customlabel{hyp1}{\Ho}]
{\it the algebra $D_{L,x}$ has no nonzero zero divisor,}

\item[\customlabel{hyp2p}{\Htpl{\pl}}\customlabelquiet{hyp2pi}{\Htpl{\pl_i}}]
{\it for all places $\ql$ above $\pl$, there exists $u_{\ql} \in 
L_{\ql}^\times$ such that
$v_{\ql}(u_\ql) = \frac{e_{\pl}} r \cdot v_{\pl}(x)$
and
$$x = \prod_{\ql | \pl} N_{L_{\ql}/K_{\pl}}(u_\ql),$$}

\item[\customlabel{hyp2}{\Ht}]
{\it for all $i \in \{1, \ldots, s\}$, the hypothesis
\ref{hyp2pi} holds.}
\end{itemize}

\bigskip

We recall that a place $\pl$ is called \emph{inert} if there is a unique place $\ql$ above $\pl$, with $e_\ql=1$.
\begin{lemma}
\label{lem:H1}
The hypothesis \ref{hyp1} holds as soon as there exists a place $\pl \in
X^\star$ which is inert in $Y$ and at which $v_\pl(x)$ is coprime 
with~$r$.
\end{lemma}

\begin{proof}
Let $\pl$ be a place satisfying the requirements of the lemma.
We embed $D_{L,x}$ into $D_{L_\pl, x} = K_\pl \otimes_K
D_{L,x}$. By Proposition~\ref{prop:divisionlocal}, we know 
that the latter has no nonzero zero divisor. The lemma follows.
\end{proof}

The hypothesis \ref{hyp2p} clearly implies that $v_{\pl}(x)$
has to be divisible by $\frac r{e_\pl}$. The next lemma shows that 
the converse is true for unramified places over a finite field.

\begin{lemma}
\label{lem:H2}
We assume that $k$ is a finite field.
Let $\pl$ be a place of $X$.
If $\pl$ is unramified in $Y$ and $v_\pl(x)$ is divisible by $r$, 
then \ref{hyp2p} holds.
\end{lemma}

\begin{proof}
Let $m_\pl$ be the number of places of $Y$ above $\pl$.
Let $\ql$ be a place over $\pl$. By assumption, the extension
$L_\ql / K_\pl$ is unramified of degree $d_\pl = r / m_\pl$.
Since, moreover, the residue field on $K_\pl$ is finite, we conclude 
that any element of $K_\pl$ of valuation divisible by $d_\pl$ is a
norm in the extension $L_\ql / K_\pl$.

Since $r$ divides $v_\pl(x)$, one can write $x$ as a product
$x = \prod_{\ql|\pl} x_\ql$ where each $x_\ql \in K_\pl$ has 
valuation $v_\pl(x) / m_\pl$. For each place $\ql$ above $\pl$,
one can then find $u_\ql \in L_\ql$ such that $N_{L_\ql/K_\pl}(u_\ql) = 
x_\ql$. This equality implies in particular that
$$v_\ql(u_\ql) = \frac {v_\pl(x_\ql)}{d_\pl} 
= \frac {v_\pl(x)} {m_\pl d_\pl} = \frac{v_\pl(x)} r.$$
On the other hand, by construction, we have
$\prod_{\ql|\pl} N_{L_\ql/K_\pl}(u_\ql) = 
\prod_{\ql|\pl} x_\ql = x$, which finally ensures that the
hypothesis \ref{hyp2p} is fulfilled.
\end{proof}

We are now ready to define our code.
For $i \in \{1, \ldots, s\}$, we consider the $k$-algebras 
$V_i \coloneqq \OO_{L_i}/t_i \OO_{L_i}$ which are finite dimensional of
dimension~$r$.
We form the $k$-vector space
$$\HH \coloneqq \End_k(V_1) \times \End_k(V_2) \times
\cdots \times \End_k(V_s).$$
which is the ambient space in which our code will eventually sit.
We equip $\HH$ with the so-called \emph{sum-rank weight} $w_{\srk}$
defined as in the Introduction by
$$w_{\srk}(\varphi_1, \ldots, \varphi_s) \coloneqq \sum_{i=1}^s
\rk (\varphi_i).$$
We now assume the hypothesis \ref{hyp2}. For each $i$, we choose
a family of elements $u_{i,\ql}$ indexed by the places $\ql$ above
$\pl_i$ satisfying the requirements of \ref{hyp2pi}. We form the
element $u_i = (u_{i,\ql})_{\ql | \pl_i} \in L_i$.
By Lemma~\ref{lem:isomDLx}, we have an isomorphism
$$\varepsilon_i :
D_{L_i,x} \stackrel{\gamma_{u_i}}\longrightarrow
D_{L_i,1} \stackrel{\varepsilon}\longrightarrow
\End_{K_i}(L_i),$$
and Lemma~\ref{lem:isomLambdaLx} indicates moreover that 
$\gamma_{u_i}$ induces an isomorphism $\Lambda_{L_i,x} \to
\Lambda_{L_i,1}$.

Take $f\in D_{L_i,x}$. Unrolling the definitions, we realize that $\varepsilon_i(f) =
f(u_i \Phi)$; hence the morphism $\varepsilon_i$ can be thought of
as the evaluation map at $u_i \Phi$.
It follows moreover from the definitions (see Subsection~\ref{ssec:DLxlocal})
that $\varepsilon_i(f)$ stabilizes the lattice $\OO_{L_i}$ whenever 
$f \in \Lambda_{L_i,x}$. For those $f$, we let 
$\bar \varepsilon_i(f) \in \End_k(V_i)$ be the reduction of 
$\varepsilon_i(f)$ modulo $t$. Noticing finally that the assumption
that $\pl_i \not\in \pi(\supp(E))$ ensures that the Riemann--Roch space $\Lambda_{L,x}(E)$ (see Definition \ref{def:RROre}) is included in $\Lambda_{L_i, x}$, we define the ``multi-evaluation'' map
\begin{equation}
\label{eq:multieval}
\begin{array}{rcl}
\alpha : \quad 
\Lambda_{L,x}(E) & \longrightarrow & \HH \\
f & \mapsto & \big(\bar \varepsilon_i(f)\big)_{1 \leq i \leq s}.
\end{array}
\end{equation}

\begin{definition}
\label{def:agsrk}
The code $\CC(x; E; \pl_1, \ldots, \pl_s)$ 
is defined as the image of $\alpha$.
\end{definition}

\begin{remark}
The code $\CC(x; E; \pl_1, \ldots, \pl_s)$ depends on the
choice of the $u_i$'s for $i=1,\dots,s$. However, this dependence is quite weak, in
the sense that changing the $u_i$'s will eventually result in a code 
which is conjugated to the initial one by an element of $\prod_{i=1}^s 
\text{GL}(V_i)$.
That is the reason why we prefer omitting to mention the $u_i$'s in
the notation.
\end{remark}

\subsection{Code's parameters}\label{ssec:parameters}

For a $k$-linear code $\CC$ sitting inside $\HH$, we define:
\begin{itemize}
\item its \emph{length} $n$ as the $k$-dimension of the ambient space
$\HH$, \ie $n \coloneqq s r^2$,
\item its \emph{dimension} $\kappa$ as its $k$-dimension, \ie 
$\kappa \coloneqq \dim_k \CC$,
\item its \emph{minimum distance} $d$ as the minimal sum-rank weight
of a nonzero codeword in~$\CC$.
\end{itemize}
Those parameters are related by the Singleton inequality which
reads $r d + \kappa \leq n + r$ in our setting~\cite[Thm.~3.1]{BGLR21}. The next theorem provides explicit lower bounds for  the dimension and the minimum 
distance of our codes.

\begin{theorem}\label{th:codeparameters}
We keep the previous notations. We assume \ref{hyp1} and \ref{hyp2}, and 
that $\deg_Y(E) < s r$.
Then, the dimension $\kappa$ and the minimum distance $d$ of
$\CC(x; E; \pl_1, \ldots, \pl_s)$ satisfy
\begin{align*}
\kappa & 
\geq r {\cdot} \deg_Y(E) - r{\cdot}(g_Y - 1)
  - \frac {r^2} 2 \sum_{\pl \in X^\star} \frac{b_\pl{-}1}{b_\pl e_\pl} \deg_X(\pl), \\
d & \geq sr - \deg_Y(E).
\end{align*}
\end{theorem}

\begin{proof}
Let $f\in \Lambda_{L,x}(E)$ be a nonzero function and set $\omega$ to be 
the sum-rank weight of $\alpha(f)$, where $\alpha$ is the evaluation 
map defined in Equation~\eqref{eq:multieval}. By definition, we have
$\sum_{i=1}^s \rk\:\bar \varepsilon_i(f) =\omega$,
where we recall that the $\bar \varepsilon_i$'s are the components of
$\alpha$. We set $d_i \coloneqq \dim_k \ker 
\bar \varepsilon_i(f)$ for $i\in \{1,\ldots,s\}$. 
By standard linear algebra, we get
\begin{equation}
\label{eq:sumdi}
\sum_{i=1}^s d_i = \sum_{i=1}^s \dim_k V_i - \rk\:\bar \varepsilon_i(f) = sr-\omega.
\end{equation}
Recall that $E = \sum_{\ql \in Y^\star} n_\ql \ql$.
We introduce the divisor
$$E' \coloneqq -\sum_{i=1}^s d_i \pl_i + 
\sum_{\pl \in X^\star} \left\lfloor\sum_{\ql|\pl} \frac{r{\cdot}n_\ql}{e_\pl m_\pl}\right\rfloor \pl
\in \Div(X),$$
where $e_\pl$ and $m_\pl$ 
were defined in 
Subsection~\ref{ssec:RROre}. It follows from Lemma~\ref{lem:isomLambdaLx}
and Propositions~\ref{prop:vtNrd1} and~\ref{prop:vtNrd2} that $\Nrd(f) 
\in L_X(E')$. Besides, we have
\begin{align*}\deg_Y(E') 
& \leq - \sum_{i=1}^s d_i + \sum_{\ql \in Y^\star} \frac{r{\cdot}n_\ql}{e_\pl m_\pl} \deg_X(\pi(\ql))\\
 &= \omega - sr + \deg_Y(E),
 \end{align*}
the last equality coming from Equation~\eqref{eq:sumdi} and the relation
$e_\pl m_\pl \deg_Y(\ql) = r {\cdot} \deg_X(\pi(\ql))$.
As a consequence, if $\omega < sr - \deg_Y(E)$, we have $\Nrd(f) = 0$.
Since $\Nrd(f)$ is, by definition, the determinant of the map
$\mu_f : D_{L,x} \to D_{L,x}$, $g \mapsto gf$, its vanishing implies
that $\mu_f$ is not injective. In other words, $f$ is a zero divisor
in $D_{L,x}$.
Thanks to hypothesis~\ref{hyp1}, we conclude that $f$ has to vanish.
In conclusion, we showed that $\omega\geq sr-\deg_Y(E)$, hence the bound 
on $d$.

As a byproduct of what precedes, we obtain the injectivity of $\alpha$.
Therefore $\kappa = \dim_k \Lambda_{L,x}(E)$, and the announced lower
bound on $\kappa$ now follows from Corollary~\ref{cor:RROre}.
\end{proof}

Putting together the bounds from the previous theorem we can highlight the \emph{defect} of our codes with respect to the Singleton bound.
\begin{corollary}\label{cor:defect}
Under the assumptions of Theorem~\ref{th:codeparameters}, and still
writing $n$, $\kappa$ and $d$ for the length, the dimension and the
minimum distance of $\CC(x; E; \pl_1, \ldots, \pl_s)$, respectively, 
we have
$$r d + \kappa \geq n + r - \left(r{\cdot}g_Y
  + \frac {r^2} 2 \sum_{\pl \in X^\star} \frac{b_\pl{-}1}{b_\pl e_\pl} \deg_X(\pl)\right).$$
\end{corollary}

\begin{remark}
\label{rk:perspective}
More generally, one can consider $k$-subspaces $W_i \subset V_i$
and replace $\alpha$ by the restricted multi-evaluation map
$$\begin{array}{rcl}
\Lambda_{L,x}(E) & \longrightarrow & 
  \Hom_k(W_1,V_1) \times \cdots \times \Hom_k(W_s,V_s) \\
f & \mapsto & \left(\bar \varepsilon_i(f)_{|W_i}\right)_{1 \leq i \leq s}.
\end{array}$$
Doing so, we obtain more general codes, for which the bounds of
Theorem~\ref{th:codeparameters} stay valid.
\end{remark}

\subsection{The case of isotrivial covers}\label{ssec:isotrivial}

Let $\ell$ be a finite cyclic extension of $k$ of order~$r$.
Given $X$ as before, we consider the curve $Y = \text{Spec }\ell 
\times_{\text{Spec }k} X$. We have a map $\pi: Y\to X$ coming from the canonical map from $\text{Spec }\ell$ to $\text{Spec }k$. The theory developed earlier applies to this setting.
In this particular case, we notice that:
\begin{enumerate}[(1)]
\item the cover $\pi : Y \to X$ is unramified everywhere, \ie 
$e_\pl = 1$ for all places $\pl \in X^\star$,
\item the Riemann--Hurwitz formula \cite[Thm.~3.4.13]{S09} asserts that
$g_Y - 1 = r \cdot (g_X - 1)$,
\item all rational places of $X$ are inert in $Y$ and, more generally,
all places whose residue field is linearly disjoint from $\ell$ are
inert; however, the reader should be careful that
if $\pl$ is an inert place of $X$ and $\ql$ is the unique place
above $\pl$, we have $\deg_Y(\ql) = r{\cdot}\deg_X(\pl)$,
\item the residue field of any place of $Y$ is a $\ell$-algebra;
in particular the codes $\CC(x; E; \pl_1, \ldots, \pl_s)$ are always
$\ell$-linear, \ie they are $\ell$-subvector spaces of $\HH$.
\end{enumerate}
In this setting, it is relevant to work with the $\ell$-length and
the $\ell$-dimension. Precisely, if $\CC$ is a $\ell$-linear code
sitting inside $\HH$, we define:
\begin{itemize}
\item its $\ell$-length $n_\ell$ as the dimension over $\ell$ of
the ambient space $\HH$, \ie $n_\ell \coloneqq sr$,
\item its $\ell$-dimension $\kappa_\ell$ by $\kappa_\ell \coloneqq
\dim_\ell \CC$.
\end{itemize}
The Singleton bound now reads $d + \kappa_\ell \leq n_\ell + 1$ where
$d$ still denotes the minimum distance.
For the code $\CC(x; E; \pl_1, \ldots, \pl_s)$ with parameters
satisfying the hypotheses~\ref{hyp1} and~\ref{hyp2}, Theorem~\ref{th:codeparameters}
provides the following lower bounds
\begin{align}
\kappa_\ell & \geq
\deg_Y(E) - r{\cdot}(g_X - 1)
  - \frac {r} 2 \sum_{\pl \in X^\star} \frac{b_\pl{-}1}{b_\pl} \deg_X(\pl), 
\label{eq:ell:dim} \\
d & \geq sr - \deg_Y(E),
\label{eq:ell:dist} 
\end{align}
from what we derive
$$d + \kappa_\ell \geq n_\ell + 1 - \left(r{\cdot}(g_X-1)+1 
  + \frac {r} 2 \sum_{\pl \in X^\star} \frac{b_\pl{-}1}{b_\pl} \deg_X(\pl)\right).$$

\begin{remark}
\label{rem:isotrivial}
We point out that, particularly in the isotrivial case, it is not 
difficult to construct a linearized Algebraic Geometry code. First, 
one chooses a curve $X$ defined over a finite field $k$, together
with a function $x$ in its function field. By choosing $x$ carefully, 
it is possible to ensure hypothesis~\ref{hyp1} and~\ref{hyp2} with
the help of the lemmas from Subsection \ref{ssec:construction}.
Precisely, one can proceed as follows. Let $h$ be the smallest 
integer not less than $2 g_X$, which is coprime with $r$; in particular this means that $2g_X\leq h \leq 2g_X+r-1$.
We take a function $x$ which has only one pole, say at $P$, of order 
$h$ exactly. This is always possible because the classical
Riemann--Roch theorem implies the existence of a function in 
$L_X(h{\cdot}P) \backslash L_X((h{-}1){\cdot}P)$. It then follows from Lemma~\ref{lem:H1}
that Hypothesis~\ref{hyp1} is satisfied given that every place is 
unramified and $v_P(x)$ is coprime with $r$.
Besides, it follows from Lemma~\ref{lem:H2} that Hypothesis~\ref{hyp2} 
will be satisfied as soon as we choose rational places $\{\pl_1,
\ldots, \pl_s\}$ at which the function $x$ does not vanish.
Since $x$ has at most $h$ distinct zeros, it is always possible to
select $s = \#X(k) - h$ such places.
Finally, one can take the divisor $E = n{\cdot} P$ (for some 
positive integer $n$) and
form the code $\CC(x; E; \pl_1, \ldots, \pl_s)$, to which the
estimations~\eqref{eq:ell:dim} and~\eqref{eq:ell:dist} apply.
\end{remark}

\subsubsection*{Linearized Reed--Solomon codes}

We now offer an explicit construction of linearized Algebraic Geometry 
codes for $X=\PP_k^1$ and $Y = \PP^1_\ell$, both viewed as curves over 
$\text{Spec }k$. We have $g_X = 0$.
We call $t$ the coordinate on $X$ and $Y$. The function fields of
$X$ and $Y$ are then $K = k(t)$ and $L = \ell(t)$, respectively.
A place $\pl\in X^\star$ (resp. $\ql\in Y^\star$) corresponds to either 
$\infty$ or to an irreducible monic polynomial in $k[t]$ (resp. in 
$\ell[t]$). A place $\pl \in X^\star$ is rational when the corresponding 
polynomial has degree $1$, \ie rational places of $X$ are in one-to-one 
correspondence with the elements in~$k$.

We choose the function $x = t \in K^\times$. For this choice, we
have $b_\pl = 1$ for all $\pl \in X^\star$, except for the places
corresponding to $0$ and $\infty$ where $b_\pl = r$. Moreover, the
algebra
$$D_{L,x} = \ell(t)[T;\Phi] / (T^r - t),$$
where $\Phi$ is a given generator of $\text{Gal}(\ell|k)$, is
canonically isomorphic to the fraction field of $\ell[T; \Phi]$.

We consider the divisor $E = \frac {m} r {\cdot} \infty\in\Div_\Q(Y)$ for a positive 
integer $m$. Coming back to the definitions, we find that the 
Riemann--Roch space $\Lambda_{L,x}(E)$ is equal to the set $\ell[T;
\Phi]_{\leq m}$ of Ore polynomials in $T$ of degree at most~$m$.

We fix rational places $\pl_1, \ldots, \pl_s$ corresponding to elements 
$c_1, \ldots, c_s \in k \sqcup \{\infty\}$. We note that they satisfy the 
hypothesis \ref{hyp2} if and only if $c_i \in N_{\ell/k}(\ell^\times)$ for all 
$i$; we assume this from now on.
The multi-evaluation morphism $\alpha$ is given by
$$\begin{array}{rcl}
\alpha : \quad 
\ell[T; \Phi]_{\leq m} & \longrightarrow & \HH \\
f & \mapsto & \big(f(u_i\Phi)\big)_{1 \leq i \leq s} ,
\end{array}$$
where $u_i \in \ell^\times$ is a preimage of $c_i$ by the norm
map. We then recover exactly the construction of linearized Reed--Solomon 
codes~\cite{MP18, CD22}. The lower bounds~\eqref{eq:ell:dim} 
and~\eqref{eq:ell:dist} specialize to $\kappa_\ell \geq m + 1$ and $d 
\geq sr - m = n_\ell - m$. The Singleton bound is then reached in this 
case, reproving that linearized Reed--Solomon codes are MSRD codes.

\subsection{Asymptotic behaviour of linearized AG codes and the Gilbert--Varshamov bound}\label{ss:GV}
In this subsection we show that our geometric codes in the sum-rank metric are asymptotically good. In particular, as in the Hamming case, we prove the existence of sequences of linearized AG codes beating the asymptotic Gilbert--Varshamov (GV) bound in the sum-rank metric.
\smallskip

In what follows we will work with codes from isotrivial covers, as developed in Subsection \ref{ssec:isotrivial}, and therefore adopt the notation from there. We will also let $q$, the power of a prime $p$, denote the cardinality of $k$. Consequently, $\ell$ is a finite field of cardinality $q^r$.

We start by giving the Gilbert--Varshamov bound as presented in \cite[Theorem 7]{OPB21}. We adapt the statement to our context, \ie to codes sitting in $\HH$, and to our notation. We recall that for a code of dimension $\kappa$, length $rs$ and minimum distance $d$, the rate is defined as $R\coloneqq \frac{\kappa}{rs}$, and the relative minimum distance as $\delta\coloneqq\frac{d}{rs}$.
\begin{theorem}[Asymptotic Gilbert--Varshamov bound]
Let $r$ be the extension degree of $\ell / k$. 
For any positive integer $s$ and any real numbers $R, \delta \in (0,1)$ with $\delta > \frac{2}{rs}$ and
\begin{align*}
R \leq \; \delta^2-\delta\left(2+\frac{2}{sr}\right)+1+\frac{2}{sr}+\frac{1}{s^2r^2}- \frac{\sum^{\delta sr-1}_{i=1} \log_q(1+\frac{s-1}{i})+\log_q(\delta sr-1)}{r^2s}-\frac{\log_q(\gamma_q)}{r^2},\end{align*}
where $\gamma_q=\prod_{i=1}^\infty (1-q^{-i})^{-1}$, there exists a sum-rank metric codes in $\End_k(\ell)^s$ of rate at least $R$ and relative minimum distance at least $\delta$.
\end{theorem}
Note that random linear codes in the sum-rank metric almost attain the GV bound with high probability \cite[Thm.~8]{OPB21}.
When we let the number of blocks $s$ go to infinity, the previous bound is asymptotically
\begin{align}\label{eq:GVinfty}
R\leq (\delta-1)^2-\frac{\delta}{r}\log_q\left(1+\frac{1}{\delta r}\right)- \frac{\log_q(1+\delta r)}{r^2}-\frac{\log_q(\gamma_q)}{r^2} + o(1),
\end{align}
where we used that
\begin{align*}
\lim_{s \to \infty} \,
\frac 1 s \sum^{\delta sr-1}_{i=1} \log_q\big({\textstyle 1+\frac{s-1}{i}}\big)
 &= \int_0^{\delta r}  \log_q\big({\textstyle 1+\frac{1}{x}}\big)\:dx \\
 &= \delta r \: \log_q\left(1+\frac{1}{\delta r}\right) + \log_q(1+\delta r).
\end{align*}

Our aim now is to show that using our linearized AG codes we can beat the previous GV bound when the length, in particular $s$, goes to infinity. To produce asymptotically good infinite sequences of our codes, we need to consider sequences of algebraic curves. 
 Let $\{X_i\}_{i\in\mathbb{N}}$ be such a sequence, $\{g_i\}_{i\in\mathbb{N}}$ be the sequence of their genera and $\{\#X_i(k)\}_{i\in\mathbb{N}}$ be the sequence of their number of rational points. 

Following the construction of Remark~\ref{rem:isotrivial}, each curve 
$X_i$
yields a linearized Algebraic Geometry code $C_i$ of length $r s_i$, dimension
$\kappa_i$ and minimum distance $d_i$ satisfying the following
inequalities:
\begin{align}
s_i & \geq \# X_i(k) - 2 g_i - r +1,
\label{eq:lengthbound} \\
d_i + \kappa_i
  & \geq r s_i - \left(r{\cdot}(g_i-1)
  + \frac {r} 2 \sum_{\pl \in X_i^\star} \frac{b_\pl{-}1}{b_\pl} \deg_{X_i}(\pl)\right).
\label{eq:ourbound}
\end{align}
We note moreover that we have complete freedom on the choice of
$\kappa_i$, which can be any integer between $0$ and $r s_i$.
By definition, $b_\pl$ is the denominator of $v_\pl(x_i)/r$ after reduction, therefore $b_\pl\leq r$. Using furthermore the fact that the
function $x_i$ that served to the construction of $C_i$ has at most $2 g_i + r$ zeros (counted with multiplicity), we obtain
  \[ \frac {r} 2 \sum_{\pl \in X_i^\star} \frac{b_\pl{-}1}{b_\pl} \deg_{X_i}(\pl)\leq \frac{(r-1)(2g_i + r-1)} 2,\]
and Equation \eqref{eq:ourbound} rewrites as
\begin{equation}\label{eq:singseq}
d_i + \kappa_i \geq r s_i - (2r{-}1) g_i - \frac{r^2 - 4r+1}2.
\end{equation}

By Equation~\eqref{eq:lengthbound}, we know that we can take $s_i=\# 
X_i(k)-2g_i-r+1$. Sending the length of our sequence of codes to 
infinity is equivalent to suppose that 
$\{\#X_i(k)-2g_i\}_{i\in\mathbb{N}}$ goes to infinity. By the 
Hasse--Weil bound for a curve $X$ of genus $g$, we have $\#X(k)-2g\leq 
q+1+2g(\sqrt{q}-1)$. Therefore, for fixed $q$, if $\#X_i(k)-2g_i$ goes to 
infinity, so must~$g_i$.

Let $A(q)$ be the Ihara's constant $A(q)$, defined by
\[A(q)\coloneqq\limsup_{g\to +\infty} \frac{\max_{\text{genus(X)=g}} \# X(k)}{g}.\]
By the work of Ihara \cite{Iha81} and Tsfasman, Vlăduţ and Zink~\cite{TVZ82}, we know that $A(q)\leq\sqrt{q}-1$. Furthermore, when $q$ is a square, Drinfeld and Vlăduţ~\cite{VD83} proved that $A(q)\geq  \sqrt{q}-1$.
Even better, still assuming that $q$ is a square, there exist sequences of algebraic curves whose number of rational points goes to infinity, and which attain the Drinfeld and Vlăduţ's bound. Examples of such sequences include modular curves (used by Tsfasman, Vlăduţ and Zink to beat for the first time the GV bound in the Hamming metric \cite{TVZ82}) and the more explicit towers of curves introduced by Garcia and Stichtenoth in~\cite{GS95}.

\begin{theorem}
\label{theo:compGV}
We assume that $q$ is a square. For all real numbers 
$R, \delta \in (0,1)$ such that
\[R < 1 -\delta - \frac{ 2}{\sqrt{q}-1}+ \frac{1}{r(\sqrt{q}-1)}\]
there exists a $\ell$-linear linearized Algebraic Geometry code
with rate at least $R$ and relative minimum distance at least $\delta$.
\end{theorem}

\begin{proof}
We consider a sequence $\{X_i\}_{i \in \N}$ attaining Drinfeld and
Vlăduţ's bound, \emph{i.e.} such that the genus of $X_i$, denoted
by $g_i$, goes to infinity and
$$\limsup_{i \to +\infty} \frac{\# X_i(k)}{g_i} = \sqrt q - 1.$$ 
For each index $i$, we set $s_i = \# X_i(k) - 2 g_i - r+1$ and $\kappa_i
= \lceil r s_i R \rceil$. By the discussion preceding the statement 
of the theorem, there exists
a code $C_i$ of length $r s_i$, dimension $\kappa_i$ and minimum
distance $d_i$ satisfying the inequality~\eqref{eq:singseq}.
We write $R_i = \frac{\kappa_i}{r s_i}$ and $\delta_i = \frac{d_i}
{r s_i}$; they are respectively the rate and the relative minimum 
distance of $C_i$. It follows from our choice of $\kappa_i$ that
$R_i \geq R$ for all $i$ and $\lim_{i \to +\infty} R_i = R$.
Moreover, dividing Equation~\eqref{eq:singseq} by $r s_i$ and 
passing to the limit, we get
$$R = \limsup_{i \to +\infty} R_i 
  \geq 1 - \liminf_{i \to +\infty} \delta_i -
  \frac{ 2}{\sqrt{q}-1}+ \frac{1}{r(\sqrt{q}-1)}.$$
Consequently, for $i$ large enough, the linearized Algebraic
Geometry code $C_i$ does satisfy the requirements of the theorem.
\end{proof}

\begin{figure}
\begin{tikzpicture}[scale=3]

\begin{scope}
\draw[-latex] (-0.1,0)--(1.15,0);
\draw[-latex] (0,-0.1)--(0,1.15);
\draw (-0.03,1)--(0.03,1);
\draw (1,-0.03)--(1,0.03);
\node[below,scale=0.8] at (1.15,0) { $R$ };
\node[left,scale=0.8] at (0,1.15) { $\delta$ };
\node[below, scale=0.7] at (1,-0.02) { $1$ };
\node[left, scale=0.7] at (-0.02,1) { $1$ };
\node[below left,scale=0.7] at (0,0) { $0$ };
\draw[BC, thick] (0, 0.833)--(0.833,0);
\draw[GV, thick] (0.000000, 1)--(0.010000, 0.983994)--(0.020000, 0.970110)--(0.030000, 0.955156)--(0.040000, 0.939602)--(0.050000, 0.923670)--(0.060000, 0.907493)--(0.070000, 0.891155)--(0.080000, 0.874718)--(0.090000, 0.858227)--(0.100000, 0.841716)--(0.110000, 0.825214)--(0.120000, 0.808743)--(0.130000, 0.792322)--(0.140000, 0.775965)--(0.150000, 0.759687)--(0.160000, 0.743499)--(0.170000, 0.727410)--(0.180000, 0.711430)--(0.190000, 0.695565)--(0.200000, 0.679823)--(0.210000, 0.664210)--(0.220000, 0.648731)--(0.230000, 0.633390)--(0.240000, 0.618193)--(0.250000, 0.603142)--(0.260000, 0.588242)--(0.270000, 0.573496)--(0.280000, 0.558907)--(0.290000, 0.544478)--(0.300000, 0.530211)--(0.310000, 0.516108)--(0.320000, 0.502172)--(0.330000, 0.488404)--(0.340000, 0.474808)--(0.350000, 0.461383)--(0.360000, 0.448132)--(0.370000, 0.435056)--(0.380000, 0.422158)--(0.390000, 0.409437)--(0.400000, 0.396895)--(0.410000, 0.384534)--(0.420000, 0.372354)--(0.430000, 0.360356)--(0.440000, 0.348542)--(0.450000, 0.336912)--(0.460000, 0.325467)--(0.470000, 0.314208)--(0.480000, 0.303135)--(0.490000, 0.292250)--(0.500000, 0.281552)--(0.510000, 0.271044)--(0.520000, 0.260724)--(0.530000, 0.250594)--(0.540000, 0.240655)--(0.550000, 0.230906)--(0.560000, 0.221348)--(0.570000, 0.211983)--(0.580000, 0.202809)--(0.590000, 0.193828)--(0.600000, 0.185040)--(0.610000, 0.176445)--(0.620000, 0.168043)--(0.630000, 0.159836)--(0.640000, 0.151823)--(0.650000, 0.144005)--(0.660000, 0.136381)--(0.670000, 0.128953)--(0.680000, 0.121721)--(0.690000, 0.114683)--(0.700000, 0.107842)--(0.710000, 0.101198)--(0.720000, 0.094749)--(0.730000, 0.088497)--(0.740000, 0.082442)--(0.750000, 0.076584)--(0.760000, 0.070923)--(0.770000, 0.065460)--(0.780000, 0.060194)--(0.790000, 0.055126)--(0.800000, 0.050256)--(0.810000, 0.045583)--(0.820000, 0.041109)--(0.830000, 0.036833)--(0.840000, 0.032756)--(0.850000, 0.028877)--(0.860000, 0.025197)--(0.870000, 0.021716)--(0.880000, 0.018433)--(0.890000, 0.015350)--(0.900000, 0.012465)--(0.910000, 0.009780)--(0.920000, 0.007294)--(0.930000, 0.005008)--(0.940000, 0.002921)--(0.950000, 0.001033);
\node[scale=0.8] at (0.6, -0.25) { $q = 7^2,\, r = 1$ };
\end{scope}

\begin{scope}[xshift=1.7cm]
\draw[-latex] (-0.1,0)--(1.15,0);
\draw[-latex] (0,-0.1)--(0,1.15);
\draw (-0.03,1)--(0.03,1);
\draw (1,-0.03)--(1,0.03);
\node[below,scale=0.8] at (1.15,0) { $R$ };
\node[left,scale=0.8] at (0,1.15) { $\delta$ };
\node[below, scale=0.7] at (1,-0.02) { $1$ };
\node[left, scale=0.7] at (-0.02,1) { $1$ };
\node[below left,scale=0.7] at (0,0) { $0$ };
\draw[BC, thick] (0, 0.85)--(0.85,0);
\draw[GV, thick] (0.000000, 1)--(0.010000, 0.980103)--(0.020000, 0.960405)--(0.030000, 0.940908)--(0.040000, 0.921610)--(0.050000, 0.902513)--(0.060000, 0.883615)--(0.070000, 0.864917)--(0.080000, 0.846420)--(0.090000, 0.828122)--(0.100000, 0.810024)--(0.110000, 0.792126)--(0.120000, 0.774428)--(0.130000, 0.756930)--(0.140000, 0.739633)--(0.150000, 0.722535)--(0.160000, 0.705637)--(0.170000, 0.688939)--(0.180000, 0.672441)--(0.190000, 0.656143)--(0.200000, 0.640045)--(0.210000, 0.624147)--(0.220000, 0.608449)--(0.230000, 0.592951)--(0.240000, 0.577653)--(0.250000, 0.562555)--(0.260000, 0.547657)--(0.270000, 0.532959)--(0.280000, 0.518461)--(0.290000, 0.504163)--(0.300000, 0.490065)--(0.310000, 0.476167)--(0.320000, 0.462469)--(0.330000, 0.448971)--(0.340000, 0.435673)--(0.350000, 0.422575)--(0.360000, 0.409677)--(0.370000, 0.396979)--(0.380000, 0.384481)--(0.390000, 0.372183)--(0.400000, 0.360084)--(0.410000, 0.348186)--(0.420000, 0.336488)--(0.430000, 0.324990)--(0.440000, 0.313692)--(0.450000, 0.302594)--(0.460000, 0.291696)--(0.470000, 0.280998)--(0.480000, 0.270500)--(0.490000, 0.260201)--(0.500000, 0.250103)--(0.510000, 0.240205)--(0.520000, 0.230507)--(0.530000, 0.221009)--(0.540000, 0.211711)--(0.550000, 0.202612)--(0.560000, 0.193714)--(0.570000, 0.185016)--(0.580000, 0.176518)--(0.590000, 0.168220)--(0.600000, 0.160122)--(0.610000, 0.152223)--(0.620000, 0.144525)--(0.630000, 0.137027)--(0.640000, 0.129729)--(0.650000, 0.122631)--(0.660000, 0.115732)--(0.670000, 0.109034)--(0.680000, 0.102536)--(0.690000, 0.096238)--(0.700000, 0.090140)--(0.710000, 0.084241)--(0.720000, 0.078543)--(0.730000, 0.073045)--(0.740000, 0.067747)--(0.750000, 0.062649)--(0.760000, 0.057750)--(0.770000, 0.053052)--(0.780000, 0.048554)--(0.790000, 0.044256)--(0.800000, 0.040157)--(0.810000, 0.036259)--(0.820000, 0.032561)--(0.830000, 0.029063)--(0.840000, 0.025764)--(0.850000, 0.022666)--(0.860000, 0.019768)--(0.870000, 0.017070)--(0.880000, 0.014571)--(0.890000, 0.012273)--(0.900000, 0.010175)--(0.910000, 0.008277)--(0.920000, 0.006578)--(0.930000, 0.005080)--(0.940000, 0.003782)--(0.950000, 0.002683)--(0.960000, 0.001785)--(0.970000, 0.001087)--(0.980000, 0.000589)--(0.990000, 0.000290);
\node[scale=0.8] at (0.6, -0.25) { $q = 11^2,\, r = 2$ };
\end{scope}

\begin{scope}[xshift=3.4cm]
\draw[-latex] (-0.1,0)--(1.15,0);
\draw[-latex] (0,-0.1)--(0,1.15);
\draw (-0.03,1)--(0.03,1);
\draw (1,-0.03)--(1,0.03);
\node[below,scale=0.8] at (1.15,0) { $R$ };
\node[left,scale=0.8] at (0,1.15) { $\delta$ };
\node[below, scale=0.7] at (1,-0.02) { $1$ };
\node[left, scale=0.7] at (-0.02,1) { $1$ };
\node[below left,scale=0.7] at (0,0) { $0$ };
\draw[BC, thick] (0, 0.8)--(0.8,0);
\draw[GV, thick] (0.000000, 1)--(0.010000, 0.980103)--(0.020000, 0.960405)--(0.030000, 0.940908)--(0.040000, 0.921610)--(0.050000, 0.902513)--(0.060000, 0.883615)--(0.070000, 0.864917)--(0.080000, 0.846420)--(0.090000, 0.828122)--(0.100000, 0.810024)--(0.110000, 0.792126)--(0.120000, 0.774428)--(0.130000, 0.756930)--(0.140000, 0.739633)--(0.150000, 0.722535)--(0.160000, 0.705637)--(0.170000, 0.688939)--(0.180000, 0.672441)--(0.190000, 0.656143)--(0.200000, 0.640045)--(0.210000, 0.624147)--(0.220000, 0.608449)--(0.230000, 0.592951)--(0.240000, 0.577653)--(0.250000, 0.562555)--(0.260000, 0.547657)--(0.270000, 0.532959)--(0.280000, 0.518461)--(0.290000, 0.504163)--(0.300000, 0.490065)--(0.310000, 0.476167)--(0.320000, 0.462469)--(0.330000, 0.448971)--(0.340000, 0.435673)--(0.350000, 0.422575)--(0.360000, 0.409677)--(0.370000, 0.396979)--(0.380000, 0.384481)--(0.390000, 0.372183)--(0.400000, 0.360084)--(0.410000, 0.348186)--(0.420000, 0.336488)--(0.430000, 0.324990)--(0.440000, 0.313692)--(0.450000, 0.302594)--(0.460000, 0.291696)--(0.470000, 0.280998)--(0.480000, 0.270500)--(0.490000, 0.260201)--(0.500000, 0.250103)--(0.510000, 0.240205)--(0.520000, 0.230507)--(0.530000, 0.221009)--(0.540000, 0.211711)--(0.550000, 0.202612)--(0.560000, 0.193714)--(0.570000, 0.185016)--(0.580000, 0.176518)--(0.590000, 0.168220)--(0.600000, 0.160122)--(0.610000, 0.152223)--(0.620000, 0.144525)--(0.630000, 0.137027)--(0.640000, 0.129729)--(0.650000, 0.122631)--(0.660000, 0.115732)--(0.670000, 0.109034)--(0.680000, 0.102536)--(0.690000, 0.096238)--(0.700000, 0.090140)--(0.710000, 0.084241)--(0.720000, 0.078543)--(0.730000, 0.073045)--(0.740000, 0.067747)--(0.750000, 0.062649)--(0.760000, 0.057750)--(0.770000, 0.053052)--(0.780000, 0.048554)--(0.790000, 0.044256)--(0.800000, 0.040157)--(0.810000, 0.036259)--(0.820000, 0.032561)--(0.830000, 0.029063)--(0.840000, 0.025764)--(0.850000, 0.022666)--(0.860000, 0.019768)--(0.870000, 0.017070)--(0.880000, 0.014571)--(0.890000, 0.012273)--(0.900000, 0.010175)--(0.910000, 0.008277)--(0.920000, 0.006578)--(0.930000, 0.005080)--(0.940000, 0.003782)--(0.950000, 0.002683)--(0.960000, 0.001785)--(0.970000, 0.001087)--(0.980000, 0.000589)--(0.990000, 0.000290);
\node[scale=0.8] at (0.6, -0.25) { $q = 11^2,\, r \to \infty$ };
\end{scope}

\begin{scope}[xshift=3.8cm, yshift=1cm, yscale=0.15, xscale=0.15]
\draw[BC, thick] (0,0)--(1,0);
\node[right,scale=0.8] at (1,0) { Theorem~\ref{theo:compGV} };
\draw[GV, thick] (0,1)--(1,1);
\node[right,scale=0.8] at (1,1) { GV bound };
\end{scope}

\end{tikzpicture}
\caption{Comparison between GV bound and Theorem~\ref{theo:compGV}}
\label{fig:GV}
\end{figure}
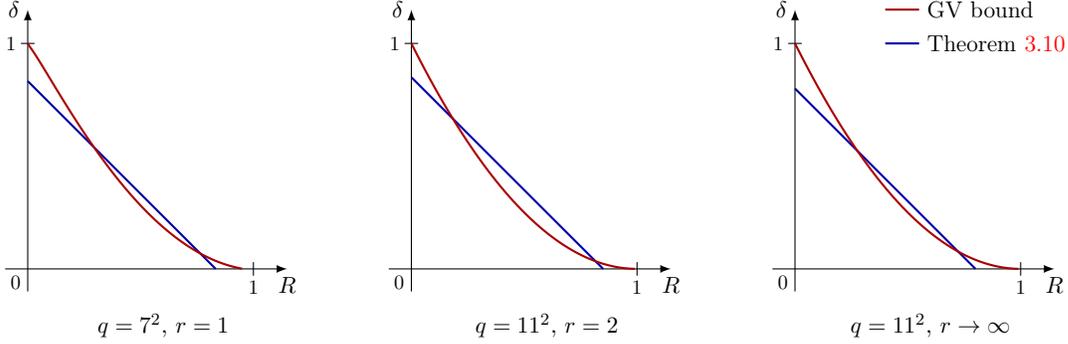

When $q$ is a square and is large enough, 
the previous bound exceeds the Gilbert--Varshamov one
for any $r$, proving 
that some families of codes in the sum-rank metric from algebraic curves 
are better than random codes. In Figure~\ref{fig:GV}, we picture the 
comparison between the GV bound and the bound from 
Theorem~\ref{theo:compGV} in several cases. We observe that,
for $q=11^2$, there is a range where our bound improves on the
Gilbert--Varshamov bound, for any $r$.
For $q = 7^2$, this only happens when $r=1$. Nevertheless,
our bounds could probably be improved by optimizing the choice of
the functions $x_i$ and the divisors~$E_i$.

\section{Conclusion}
\label{sec:conclusion}

In this article, we introduced a new family of codes for the
sum-rank metric and provided lower bounds on their dimension and
minimum distance, showing that our codes exhibit quite good
parameters. Our construction is based on algebraic geometry and
can be considered as an extension of that of AG codes to a
noncommutative framework.

\subsection{Comparison with Morandi and Sethuraman's codes}

In~\cite{MS98}, Morandi and Sethuraman proposed a construction quite similar
 to ours, whose initial input is a maximal order in a central simple algebra over the function field of a curve.
Our approach meets Morandi and Sethuraman's one because the rings $D_{L,x}$ 
we considered in this paper turns out to be central simple algebras
over $K = k(X)$ (where we recall that $X$ is a curve defined over~$k$). 
Additionally, the main ingredient in Morandi and Sethuraman's article
is a noncommutative version of the Riemann--Roch's 
theorem~\cite[Thm.~4]{MS98} (which is initially due to Van 
Geel~\cite{Geel81a, Geel81}), which looks similar to our 
Corollary~\ref{cor:RROre}. Nevertheless, our contribution differs from~\cite{MS98} in several 
important points.

First, we are working with the sum-rank distance while 
Morandi and Sethuraman work with the classical Hamming metric. Beyond 
this obvious separation, the setup of~\cite{MS98} forces the authors to 
choose ``evaluation points'' which are totally ramified places of the 
central simple algebra. In comparison, we have more freedom in our 
framework, being only constrained by the hypothesis~\ref{hyp2p}, which 
is of different nature but usually much weaker.

Secondly, Morandi and Sethuraman's construction uses \emph{maximal}
orders in the underlying division algebra whereas the explicit rings
$\Lambda_{L,x}$ we introduced are usually \emph{non-maximal} orders.
Here, the divergence is more subtle but, in some sense, it is the
same as the difference between smooth and singular curves in the
classical Riemann--Roch's theorem. 
Indeed, in the commutative setting, the ring of functions on smooth 
curves which are regular outside one fixed place is a Dedekind domain 
which is a maximal order in the corresponding field of functions. On 
the contrary, when the curve is singular, the order defined by the 
ring of regular functions is no
longer maximal (and desingularizing the curve consists in replacing
this order by the maximal one).
Following this analogy, our Corollary~\ref{cor:RROre} can be thought
as a (weak) instance of an hypothetic extension of the Riemann--Roch's
theorem for central simple algebras to the singular case.

Lastly, on the practical side, our construction looks better suited for 
concrete implementation. Indeed, the main ingredients we are using are 
Ore polynomials and Riemann--Roch spaces. Both of them are available in 
standard softwares of Symbolic Computation 
(\emph{e.g.~}\textsc{Magma}~\cite{magma}, \textsc{SageMath}~\cite{sage}),
making rather concrete the perspective of implementing our codes and potentially use them.
Instead, Morandi and Sethuraman use abstract central simple 
algebras (encoded by their Hasse invariants) and the general theory of 
maximal orders inside them. Although these objects are very important in 
algebraic geometry, as far as we know, a full support for manipulating 
them on computers is not yet available. For sake of completeness, we mention however that an 
implementation in the framework of number fields is available in 
PARI/GP~\cite{pari}, after the work of Aurel Page.

\subsection{Perspectives}
To start, we plan to understand the interactions between
our construction and that of~\cite{MS98}, and possibly to set up a
general framework which encompasses both approaches.
This is a long-term project since, as discussed above, it will require
at least to extend Van Geel's noncommutative version of Riemann--Roch's theorem to
the ``singular case''. It will also need to allow for more general
divisors. Indeed, the shape of divisors $\sum_{\ql} n_\ql \ql$ with
which we have worked in the present article looks more general than the 
divisors considered in~\cite{Geel81a, Geel81, MS98}, which were
restricted to the form $\sum_{\pl} n_\pl \pl$ (the sum is taken over
$\pl$, not over~$\ql$).

Apart from this, we plan to study the decoding problem, at least
in the case of unique decoding.
Indeed, efficient decoding algorithms are available for both AG codes and linearized Reed--Solomon codes. It is then a natural
question to try to extend those algorithms to the setting of the
present paper.

Finally, it would be desirable to have a duality theorem for the codes
$\CC(x; E; \pl_1, \ldots, \pl_s)$, in the spirit of the main result
of~\cite{CD22}. Again, this is not immediate as it will require
to develop the theory of differential forms and residues in the
framework of central simple algebras. We nevertheless plan to go
back on this question in a forthcoming article.

\section*{Acknowledgment} The authors thank the anonymous second reviewer for some of the remarks that led to the writing of Section \ref{ss:GV}. This work was funded in part by the grants ANR-21-CE39-0009-BARRACUDA and ANR-18-CE40-0026-01-CLap--CLap from the French National Research Agency. Till the end of April 2023 the first author has also received funding from the  European Union’s Horizon 2020 research and innovation programme under the Marie Skłodowska--Curie grant agreement No 899987. 

\bibliography{biblio_sumrank}
\bibliographystyle{IEEEtranS}

\end{document}